\documentclass[10pt,a4paper]{article} 

\usepackage{anysize}
\usepackage[format = hang]{caption}
\usepackage{amsmath,amsthm,verbatim,amssymb,amsfonts,amscd, graphicx}
\usepackage{graphics,natbib}
\usepackage{titlesec,footmisc}
\usepackage{hyperref} 

\usepackage{listings}
\usepackage{booktabs} 

\theoremstyle{plain}

\theoremstyle{definition}

\def\dd{{\rm d}}

\def\N{{\rm N}}

\def\eps{\varepsilon}

\def\beq{\begin{eqnarray}}
\def\eeq{\end{eqnarray}}

\def\beqn{\begin{eqnarray*}}  
\def\eeqn{\end{eqnarray*}}

\def\dd{{\rm d}}
\def\N{{\rm N}}

\def\eps{\varepsilon}

\titleformat{\section}{\normalfont\large\sc\centering}{\thesection}{1em}{}
\titleformat{\subsection}[runin]{\normalfont\large\bfseries}{\thesubsection}{1em}{}
\numberwithin{equation}{section} 
\renewenvironment{abstract}
               {\list{}{\rightmargin\leftmargin}%
                \item[\text{\hspace{10mm}\sc Abstract.}]\relax}
               {\endlist}

\begin{document}

\def\heute{March 2026}

\begingroup
\begin{centering} 

  \Large{\bf But some are more equal than others}
  \\[0.8em]
\large{\bf Nils Lid Hjort} \\[0.3em] 
\small {\sc Department of Mathematics, University of Oslo} \\[0.3em]
\small {\sc January 2017, but arXiv version
  March 2026\footnote{A FocuStat Blog post from January 2017; 
  this modified form March 2026 for arXiv and other channels}}\par
\end{centering}
\endgroup


\begin{abstract}
\small{
  All men are created equal, proclaimed Jefferson in 1776 -- but
  some are more equal than others, added Orwell in Animal Farm in 1945.
  So what's the probability that two skaters are exactly equal,
  to the third decimal places, after four distances?}
    

\end{abstract}


\section*{Are Allan and Odin exactly equal?}
\label{section:maxima}

\begin{figure}[h]
\centering
\includegraphics[scale=0.33]{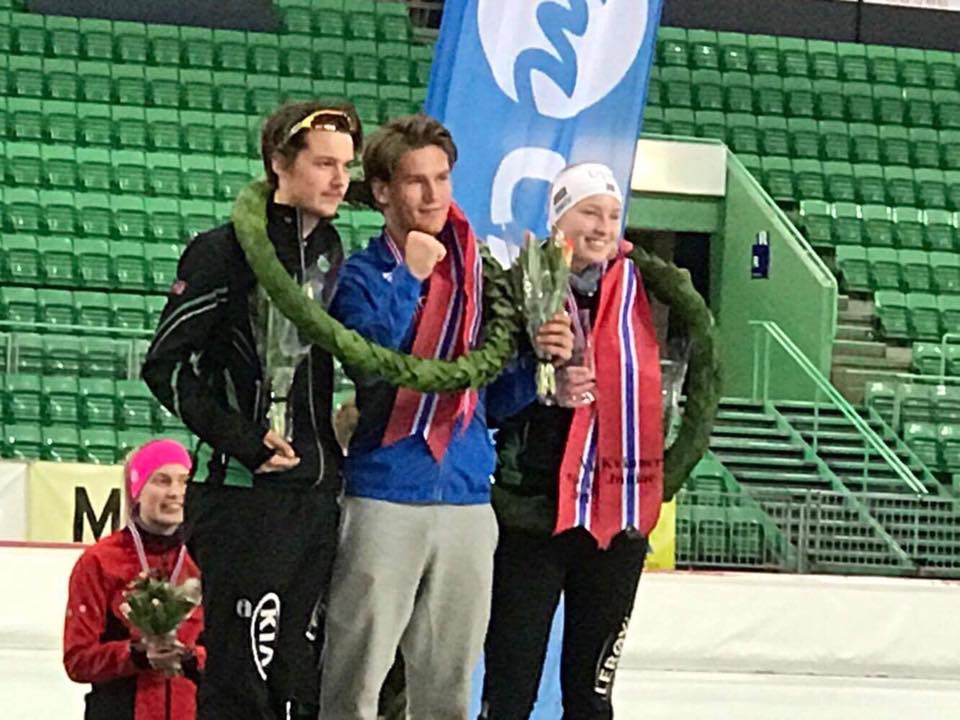}
\caption{Never before:
  shared gold, pointsums after four distances:
  Allan Dahl Johansson and Odin By Farstad,
  at the Norwegian sprint speedskating Junior Championships. 
  Ragne Wiklund, the Norwegian champion for junior women, to the right.}
\label{figure:queen1}
\end{figure}

\noindent 
The world of sports is famously able and zealously eager
to rank its competitors, even if it takes millimetres
and milliseconds and highlevel photography to sort
one medal from another. But once in a very rare while
even the ministries of sport can't tell the difference.
What's the probability that Allan and Odin are exactly
equal, to three decimal places in pointsum, after four distances?

The photo above is historic \& epic, as it reports on
an event which never has happened before, in the
World History of Speedskating: The two top skaters,
after having skated four separate distances, are equal
in pointsum, to three decimal places, and are hence
sharing the laurels. Just how wondrously improbable
\& freakish is this?

Let's tell the story first, and tend to some probability
calculus afterwards. At the 2017 Norwegian Junior Sprints
Championships (January 14-15, Vikingskipet Olympic Hall Hamar),
the two favourites were Allan Dahl Johansson and Odin By Farstad.
They are hopefully future Olympians, and the nation's eyes
(well, those of the sufficiently educated strata)
are upon them.

On the Saturday, Allan did 37.49 (3rd) and 1:12.82 (1st)
in the 500-m and 1000-m, whereas Odin managed 36.75 (1st)
and 1:13.92 (3rd). Samalogue pointsums are computed by
converting times to 500-m scale and adding, so Allan had
37.49 + 72.82/2 = 73.900 and Odin 36.75 + 73.92/2 = 73.710
by then; thus Allan needed to win back 0.19 seconds
on the Sunday. On the Sunday 500-m, Allan and Odin did 37.12 (2nd)
and 36.95 (1st), yielding pointsums 111.020 and 110.660.
For the two athletes, and for everyone present
(since the speaker is giving them the details),
it's clear that Allan needs to beat Odin with hum-ti-dum
(111.020 - 110.660)*2 = 0.72 seconds to become equal.
First comes Allan, with a clean, strong race and a personal
best of 1:12.35 (which incidentally is 0.08 seconds faster
than when Dan Jansen apparently brought
The United States of America to tears at the same rink
at the 1994 Olympics). Odin then knows he needs to set
a personal best-ever mark too, of, let's see, 1:13.07,
to equalise Allan. And so he does!,
exactly dot-on-the-spot, 1:13.07.

\begin{figure}[h]
\centering
\includegraphics[scale=0.44]{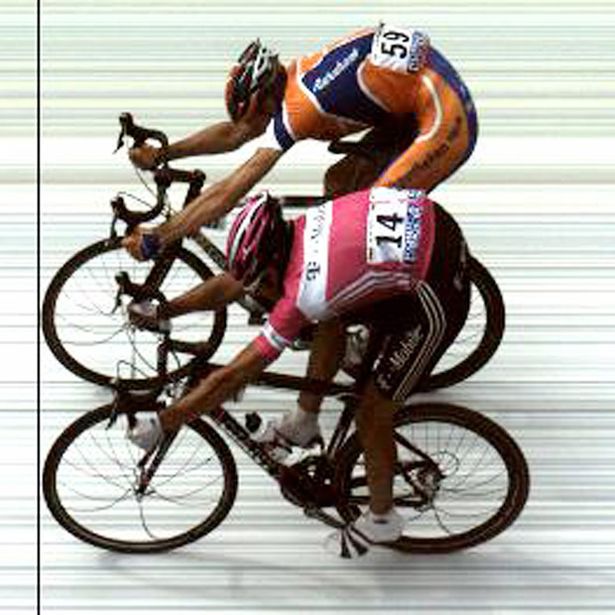}
\caption{No human eye can spot any difference when they finish.}
\label{figure:queen2}
\end{figure}

Almost by definition this is if not incredible, then
a totally freakish event, as it has never happened before.
There have been various close shaves and occasions where
only the engineering equipment manages to `see' the difference.
At the 2014 Sochi Olympics, both Zbigniew Bródka and
Koen Verweij did 1:45.00 in the 1500-m, beating all others.
In speedskating, all official times are in hundredths
of seconds (i.e.~to two decimal places, not three).
Most fans, and indeed most human beings, would've been
happy to have the Pole and the Dutch share the gold medal.
After all, when Thomas Wassberg skied the 15 km cross-country
in 41 minutes and 57.63 seconds, at the 1980 Olympics,
officials and fans quickly agreed, a posteriori,
that it was pretty ridiculous to say that Juha Mieto,
who did 41 minutes and 57.64 seconds, only deserved
the silver, so the rules were changed -- hundredths of
a second don't count anymore, in cross-country skiing.
The International Skating Union (ISU), however,
is severely and devotedly hunting for subtle differences,
even when not noticeable for the human eye. So they
went into their computers and found that Bródka did 1:45.006,
Verweij merely 1:45.009, thus splitting the horse's hair
and saying gold to Poland, silver to the Netherlands.

Luckily (I dare to express), the Hamar time-keeping technicians
didn't find any such 0.001-level differences on this occasion,
so lo \& behold, Allan and Odin were judged to be perfectly
equal. So they're forever sharing their gold \& their laurels.

\begin{figure}[h]
\centering
\includegraphics[scale=0.33]{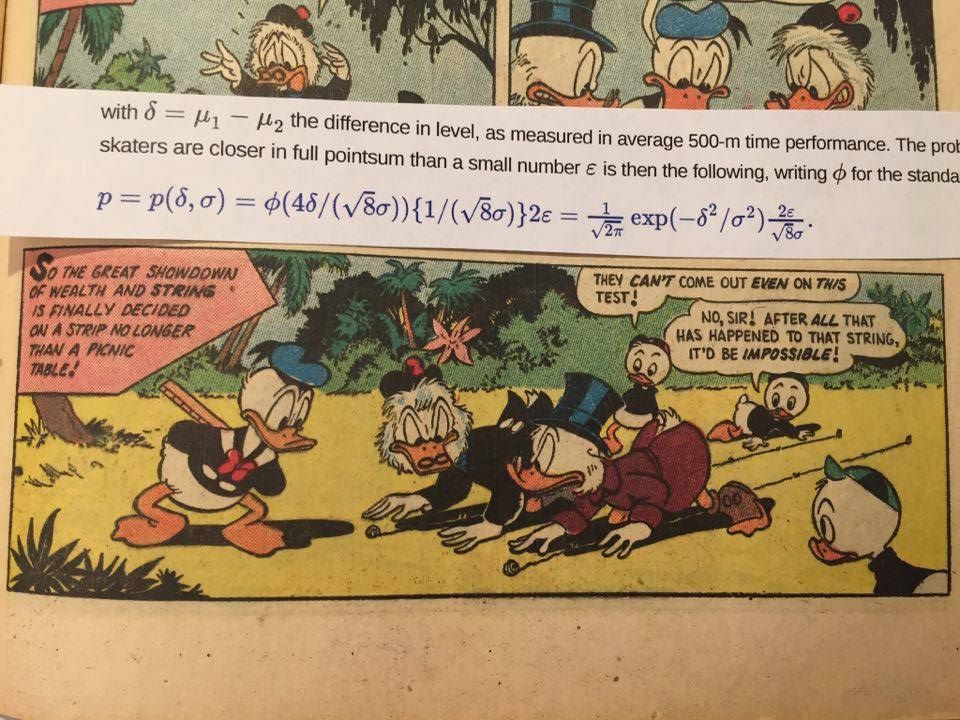}
\caption{Scrooge McDuck and Flintheart Glomgold
  fight it out, very evenly, after dramatic days
  and nights in the African wilderness;
  Carl Barks, Scrooge \#15, 1956.}
\label{figure:queen3}
\end{figure}

\section*{Close Encounters of the Statistical Kind} 

Any decent sport has a catalogue of Close Encounters.
In the last pair of the World Sprint Championship
in Trondheim, 1984, Gaétan Boucher had to beat Sergey Khlebnikov
with exactly 0.48 seconds or more, to become World Champion,
and for good measure the Canadian won over the Russian with
0.49 seconds. Even more thrillingly, Marius Bratli in the
last pair of the Norwegian Junior Championships 2016 knew
he had to skate the 5000-m in precisely 7:21.33 to beat
comrade Vetle Stangeland and become the champion -- which
is very precisely what he did, to the 0.01 of a dot.
You can't ever do this again even if you try very hard.

\begin{verbatim}
 9	An Liu       China  35.94
10	An Liu       China  35.95
11	Tao Yang     China  35.96
12	Xuefeng Sun  China  35.96
\end{verbatim}

Other freakish-looking incidents take place too.
In the Berlin World Cup, January 29, 2017, there were
four Chinese taking part in the 500-m Group B. They skated
almost exactly at the same speed (and two of them had
identical names, adding to the confusion and fascination).

Our own lives carry surprises too, and I suppose you,
gentle reader, have stories matching this one, which now
comes to my mind: Once upon a time I was travelling in the
northern and almost empty part of Sweden, on my way to
give lectures on pattern recognition at the Ammarnäs
Winter School. We stopped at a boring café with
luke-warm coffee in The Middle of Nowhere (and it was
minus 27 Celsius). At the neighbouring table was a little boy,
dutifully reading his comics magazine, with his grandmother.
I heard him say (well, I'm translating from Swedish,
for convenience), ``Grandma, have you ever heard about Professor
Nils Lid Hjort?''. She answered ``no'', incidentally,
to my modest disappointment.

Attempts at making probability calculations for such
`strange events' are difficult, in particular because
the questions are so often raised after the fact,
a posteriori. The Law of Very Large Numbers,
as investigated and formulated by famous statisticians
Persi Diaconis and Fred Mosteller in a JASA 1989 paper,
implies that ``with a large enough sample, any outrageous thing
is apt to happen'', with or without a relevant and contextual
semi-explanation, so to speak (my pattern recognition trip
to Ammarnäs in 2001 took place the only time in history
I had a three-page paper published in the Swedish Kalle Anka).
So it's clear that some 1-in-a-billion sounding chances
might in reality be less freakoutish when properly
contextualised and decoded; statisticians and probabilists
bring in their `what are we conditioning on?' themes.

Assessing the rareness of the Allan = Odin type incidences
is however an approachable problem. For two top-level
and almost equally talented skaters, let's look at their
final pointsums after four distances, say $X=X_1+X_2+X_3+X_4$
and $Y=Y_1+Y_2+Y_3+Y_4$.
A good approximation to skating reality is that the $X_j$
and $Y_j$ come from normal distributions, with parameters say
$(\mu_1,\sigma^2)$ and $(\mu_2,\sigma^2)$,
and that these are independent. Hence, for the difference,
\beqn
Z=X-Y\sim\N(4\delta,8\sigma^2), 
\eeqn 
with $\delta=\mu_1-\mu_2$ the difference in level,
as measured in average 500-m time performance.
The probability that the two skaters are closer in full
pointsum than a small number $\eps$
is then the following, writing $\phi$ for the standard normal density:
\beqn
p=p(\delta,\sigma)=\phi\Bigl({4\delta\over \sqrt{8}\sigma}\Bigr)
   {1\over \sqrt{8}\sigma}2\eps
   ={1\over \sqrt{2\pi}}\exp\Bigl(-{\delta^2\over \sigma^2}\Bigr)
   {2\eps\over \sqrt{8}\sigma}. 
\eeqn 

Here we first put in $\eps=0.005$, since pointsums closer
than this means equality, for the 500-1000-500-1000 combination
with time results on the 0.01 scale. A decent value for
$\sigma$, meant to reflect the standard deviation level
in repeated top races by the same skater, again on the 500-m time
scale, is 0.50 (an exercise left to the reader is to estimate
this standard deviation based on averaging the empirical variances
for Allan and Odin, for the eight races noted above,
and the answer is 0.460). If we now further assume
that the two skaters are in principle of the very same
top level of ability, then the resulting probability is 0.00282,
or 2.82 per mille. If Allan and Odin should be willing to do this,
over and over again, over the coming 350 weekends or so
(it'll take some seven years, and we can watch),
we can expect that they are this close, close enough
to share the gold medal, in 1 of these 350 competitions.

If the two top skaters we're comparing are not of the
mathematically identical level of ability, a different formula
ensues, and the probability in question becomes smaller.
If we take $\delta\sim \N(0,\tau^2)$, with $\tau$ 
indicating the degree to which the two skaters can be
expected to be different, the probability of equal
pointsums in the end becomes
\beqn
p=\int p(\delta,\sigma)g(\delta)\,\dd\delta
={1\over \sqrt{2\pi}}{2\eps\over \sqrt{8}\sigma}
   {1\over \sqrt{1+2\tau^2/\sigma^2}}. 
\eeqn 

For the situation with two top skaters of very nearly
but not identical ability levels, this might indicate that
the Allan-Odin probability is a bit lower, perhaps 2.30
per mille (reflecting a $\tau$ of around a quarter of a second),
rather than 2.82 per mille.

During the regular and always far-ranging discussions at
{\it Forum for Skøytehistorie}
(``ein av dei mest begeistra, intense og nerdete stader
som finnest på nettet'', according to Dagbladet)
it was mentioned that Allan's and Odin's pointsums
were equal even down to the 0.001 level, following computer
clock checks of the result times (these are not official).
In that case we may put $\eps=0.001$ above, 
which leads to the even more impressively slim probability
of 0.56 per mille, for the case of two perfectly matched
skaters, and 0.43 per mille for the case of not perfectly
even skaters, with $\tau=0.25$ above.

As perhaps illustrated here, statisticians have a certain
fascination for interesting rare events -- which adds
one more reason for us to hope that Allan and Odin
share the 1500 m gold medal in the 2026 Olympics in Milano.
After all, strange things are bound to happen -- see also
our blogposts about the two Lillestrøm women born
on the same day who died on the same day, 101 years later,
and about conspiracy probability calculus for footballers.

\section*{References} 

\parindent0pt
\parskip3pt

C Barks (1956). 
The Second-Richest Duck.
{\it Uncle Scroge} \#15, Disney;
published in 26 countries; 
republished in the US 14 times; 
Figure 2 is from the original 1956 edition. 

P Diaconis, F Mosteller (1989). 
Methods for Studying Coincidences.
{\it Journal of the American Statistical Association},
84, 853--861. 

M Jullum, DR Baños, NL Hjort (2015). 
To liv: kvinnene i Lillestrøm som ble født på samme dag
og døde på samme dag.
{\it FocuStat Blogpost}, Department of Mathematics, University of Oslo.

NL Hjort (2015). 
Conspiracy Probability Calculus for Norwegian Footballers.
{\it FocuStat Blogpost}, Department of Mathematics, University of Oslo.



\end{document}